\documentclass[11pt]{amsart}
\usepackage{amsfonts,amssymb}
\usepackage[mathscr]{eucal}
\usepackage{amsmath, amsthm}
\usepackage{mathrsfs}
\usepackage{amsbsy}
\usepackage{dsfont}
\usepackage{bbm}
\makeindex
\input xypic
\xyoption{all}
\begin{document}
\baselineskip = 5mm
\newcommand \lra {\longrightarrow}
\newcommand \hra {\hookrightarrow}
\newcommand \ZZ {{\mathbb Z}} 
\newcommand \NN {{\mathbb N}} 
\newcommand \QQ {{\mathbb Q}} 
\newcommand \RR {{\mathbb R}} 
\newcommand \CC {{\mathbb C}} 
\newcommand \sgr {{\mathfrak S}} 
\newcommand \bcA {{\mathscr A}}
\newcommand \bcB {{\mathscr B}}
\newcommand \bcC {{\mathscr C}}
\newcommand \bcD {{\mathscr D}}
\newcommand \bcE {{\mathscr E}}
\newcommand \bcF {{\mathscr F}}
\newcommand \C {{\mathscr C}}
\newcommand \X {{\mathscr X}}
\newcommand \im {{\rm im}}
\newcommand \Hom {{\rm Hom}}
\newcommand \colim {{{\rm colim}\, }} 
\newcommand \iHom {{\underline {\rm Hom}}} 
\newcommand \End {{\rm {End}}}
\newcommand \Ob {{\rm {Ob}}}
\newcommand \coker {{\rm {coker}}}
\newcommand \id {{\rm {1}}}
\newcommand \DM {{\mathscr {DM}}} 
\newcommand \cm {{\sf CHM}}
\newcommand \cma {{\sf CHM_{\leq 1}^{\otimes }}}
\newcommand \ccor {{\mathscr {Corr}}}
\newcommand \uno {{\mathbbm 1}} 
\newcommand \Le {{\mathbbm L}} 
\newcommand \trf {{\rm {tr}}} 
\newcommand \AF {{\mathbb A}} 
\newcommand \PR {{\mathbb P}} 
\newcommand \Spec {{\rm {Spec}}}
\newcommand \Sm {{\mathscr {Sm}}} 
\newcommand \SP {{\mathscr {S}\! \! \mathscr {P}}} 
\newcommand \Sch {{\mathscr {Sch}}} 
\newcommand \Pic {{\rm {Pic}}}
\newcommand \FPic {\underline {{\rm {Pic}}}}
\newcommand \FDiv {\underline {{\rm {Div}}}}
\newcommand \ch {{CH}}
\newcommand \Corr {{Corr}}
\newcommand \Sym {{\rm {Sym}}}
\newcommand \spd {{\rm {s}}}
\newcommand \cha {{\rm {char}}}
\newcommand \tr {{\rm {tr}}} 
\newcommand \res {{\rm {res}}} 
\newtheorem{theorem}{Theorem}
\newtheorem{lemma}[theorem]{Lemma}
\newtheorem{sublemma}[theorem]{Sublemma}
\newtheorem{corollary}[theorem]{Corollary}
\newtheorem{example}[theorem]{Example}
\newtheorem{exercise}[theorem]{Exersize}
\newtheorem{proposition}[theorem]{Proposition}
\newtheorem{remark}[theorem]{Remark}
\newtheorem{notation}[theorem]{Notation}
\newtheorem{definition}[theorem]{Definition}
\newtheorem{conjecture}[theorem]{Conjecture}
\newenvironment{pf}{\par\noindent{\em Proof}.}{\hfill\framebox(6,6)
\par\medskip}
\title{\bf Motives of smooth families and cycles on threefolds}
\author{V. Guletski\v \i
}

\date{December 13, 2005}

\begin{abstract}
\noindent Let $X\to S$ be a smooth projective family of surfaces
over a smooth curve $S$ whose generic fiber $X_{\eta }$ is a surface
with $H^2_{et}(X_{\bar \eta },\QQ _l(1))$ spanned by divisors on
$X_{\eta }$ and $H^1_{et}(X_{\bar \eta },\QQ _l)=0$. We prove that,
if the motive of $X/S$ is finite dimensional, the Chow group $\ch
^2(X)_{\QQ }$ is generated by a multisection and vertical cycles,
i.e. one-dimensional cycles lying in fibers of the above map. If
$S=\PR ^1$, then $\ch ^2(X)_{\QQ }=\QQ \oplus \QQ ^{\oplus n}$,
where $n\leq b_2$ and $b_2$ is the second Betti number of the
generic fiber. Generators in $CH^2(X)$ can be concretely expressed
in terms of spreads of algebraic generators of $H^2_{et}(X_{\bar
\eta })$. We also show where such families are naturally arising
from, and point out the connection of the result with Bloch's
conjecture.
\end{abstract}



\maketitle

\section{Introduction}
\label{intro}

Let $X$ be a smooth projective variety defined over a field. An
important problem in algebraic geometry is to compute the Chow group
$CH^j(X)$ of codimension $j$ algebraic cycles modulo rational
equivalence relation on $X$. In codimension one the problem is
solved in the following sense: algebraically trivial divisors are
parametrized by the Picard variety of $X$. But the next step, when
$j=2$, leads to hard problems. Even if $X$ is a surface,
zero-dimensional cycles on $X$ are not understood by know. The
conjecture due to S.Bloch asserts that, if $X$ is a complex surface
with no non-trivial globally holomorphic $2$-forms, the kernel of
the Albanese mapping is trivial, \cite{Jannsen2}. If $X$ is of
special type, i.e. the Kodaira dimension of $X$ is less than two,
the conjecture was done in \cite{BKL}. When $X$ is of general type,
the problem is opened, except for a few cases (see \cite{Barlow},
\cite{InoseMizukami} and \cite{Voisin}).

In the last years it was discovered that Bloch's conjecture is
closely connected with the tensor structures in motivic categories.
Namely, in \cite{GP2} we have shown that the Albanese kernel is
trivial for a given complex surface $X$ with $p_g=0$ if and only if
its motive $M(X)$ is finite dimensional in the sense of S.Kimura,
\cite{Kimura}. The use of Kimura's theory arises, actually, from
simple facts in representation theory of symmetric groups applied in
the setting of tensor categories. A purely algebraic version of that
theory was independently developed by P.O'Sullivan, \cite{AK}. The
connection with geometry is provided by the following nilpotency
theorem due to S.Kimura, \cite[7.5]{Kimura}: any numerically trivial
endomorphism of a finite dimensional Chow motive is nilpotent.

The goal of the present paper is to show that motivic finite
dimensionality can be also useful in the study of codimension two
cycles on threefolds. In fact, we will try to extend the arguments
from \cite{GP2} to one-parameter families of smooth projective
surfaces whose generic fiber has an algebraic second Weil cohomology
group. For that we will use a more general version of Kimura's
theorem proved by Andre and Kahn, \cite[9.1.14]{AK}: any numerically
trivial endomorphism of a finite dimensional object in a nice tensor
category is nilpotent.

To state the result we need to fix some notation. All Chow groups
will be with coefficients in $\QQ $. Let $k$ be a field of
characteristic zero, and let
   $$
   \gamma :X\lra S
   $$
be a smooth projective family of surfaces over a smooth connected
quasi-projective curve $S$ over $k$. Define the Chow group
  $$
  CH^2(X)_0=\ker (\gamma _*:CH^2(X)\lra CH_1(S)=\QQ )
  $$
of $1$-cycles of degree zero with respect to the map $\gamma $, so
that
  $$
  CH^2(X)=CH^2(X)_0\oplus \QQ \; .
  $$
A one-dimensional cycle class on $X$ is called vertical, if it can
be represented by a linear combination of curves lying in closed
fibers of the map $\gamma $. Let $H^*(-)$ be a Weil cohomology
theory over the function field $k(S)$, say $l$-adic \'etale
cohomology groups. In particular, if $X_{\eta }$ is the generic
fiber of the map $\gamma $, then
  $$
  H^*(X_{\eta })=H^*_{et}(X_{\bar \eta },\QQ _l)\; ,
  $$
where $\bar \eta $ is the spectrum of an algebraic closure of
$k(S)$. Let
  $$
  H^2_{\tr }(X_{\eta })
  $$
be the transcendental part of $H^2(X_{\eta })$, i.e. the second
direct summand of $H^2(X_{\eta })$ after splitting of classes of
divisors on $X_{\eta }$. If $H^2_{\tr }(X_{\eta })=0$, it follows
that $H^2(X_{\eta })$ is generated by divisors $D_1,\dots ,D_{b_2}$
on $X_{\eta }$, where $b_2$ is the second Betti number of the
generic fiber. For each $i$ let
  $$
  W_i
  $$
be a spread of $D_i$ over $S$ (see \cite[4.2]{GreenGriffiths} or
Section \ref{spreads} below for discussion of spreads). Our main
result is then as follows:

\begin{theorem}
\label{th1} Let $\gamma :X\to S$ be a family as above, and assume
that
  $$
  H^1(X_{\eta })=H^2_{\tr }(X_{\eta })=0\; .
  $$
Then, if the relative motive $M(X/S)$ is finite dimensional, any
cycle class in $CH^2(X)_0$ is vertical, and it can be represented by
a linear combination of intersections of the spreads $W_i$ with
closed fibers of the map $\gamma $. If, moreover, $S=\PR ^1$, then
  $$
  CH^2(X)_0=\QQ ^{\oplus n}
  $$
with
  $$
  n\leq b_2\; ,
  $$
and any degree zero one-dimensional algebraic cycle on $X$ is
rationally equivalent to a linear combination of cycles
  $$
  W_1\cdot F,\; \dots ,W_{b_2}\; \cdot F\; ,
  $$
where $F$ is a fixed closed fiber of $\gamma $.
\end{theorem}

Geometrically this result can be illustrated as follows. Assume
$S=\PR ^1$ and suppose we are given, for example, two curves $C_1$
and $C_2$ on $X$ projecting onto $\PR ^1$ with the same degree. Then
the difference $C_1-C_2$ is rationally equivalent to a linear
combination of one-cycles $W_1\cdot F,\dots ,W_{b_2}\cdot F$ on $X$.

Bloch's conjecture may be considered as a problem both on
zero-dimensional cycles, and on codimension two cycles. If we look
on it from the second viewpoint, a three-dimensional counterpart of
Bloch's conjecture, in the sense of a smooth one-parameter family of
projective surfaces $X\to S$, may be stated as follows: if $H^2_{\tr
}(X_{\eta })=0$ and $H^1(X_{\eta })=0$, then $CH^2(X)_0$ is
generated by one-dimensional cycles lying in fibers of the map $X\to
S$. If we ``collapse" the curve $S$ into a point, the vertical
generators disappear, so that we obtain the usual 2-dimensional
situation. Thus, Theorem \ref{th1} partially supports the
three-dimensional conjecture in those cases where we can prove
finite dimensionality of $M(X/S)$. But we cannot avoid the
restriction $H^1(X_{\eta })=0$ yet, because we do not know any
appropriate construction of the relative Murre decomposition for a
smooth projective family of surfaces with non-trivial irregularity.
On the contrary, if the first cogomology of the generic fiber
vanishes, the relative Murre decomposition for the whole family can
be easily constructed adding some vertical correcting terms to
standard projectors. This method allows also to generalize Kimura's
theorem, \cite[4.2]{Kimura}, on motives of smooth projective curves
over a field (Proposition \ref{relcurves}).

The paper is organized as follows. In Section \ref{prelim} we recall
some well known facts about Chow motives, finite dimensional objects
and spreads of algebraic cycles. Section \ref{unole} is devoted to
splittings of the unit and a tensor power of the Lefschetz motive
from the Chow motive of a family of relative dimension $e$. The
central Section \ref{surfaces} contains the proof of Theorem
\ref{th1}. In Section \ref{curves} we generalize Kimura's theorem on
motives of curves. In the last Section \ref{examples} we show how to
construct families satisfying the assumptions of Theorem \ref{th1}
by means of spreadings out of surfaces defined over $\CC $.

\medskip

{\sc Acknowledgements.} The author wishes to thank I. Panin for a
useful remark originated this paper. The work was supported by the
Von Neumann and the Ellentuck Funds, and by NSF Grant DMS-0111298. I
am grateful to the Institute for Advanced Study at Princeton for the
support and excellent hospitality in 2005, and to JSPS, Tokyo
University and Hiroshima University for the support of my stay in
Japan, June - July 2005.

\section{Preliminary results}
\label{prelim}

\subsection{\it Chow motives over a base}

Let $S$ be a smooth connected quasi-projective variety over a field
$k$, and let $\SP (S)$ be the category of all smooth and projective
schemes over $S$. Assume we are given two objects of that category,
and let $X=\cup _jX_j$ be the connected components of $X$. For any
non-negative $m$ let
  $$
  \Corr _S^m(X,Y)=\oplus _jCH^{e_j+m}(X_j\times _SY)
  $$
be the group of relative correspondences of degree $m$ from $X$ to
$Y$ over $S$, where $e_j$ is the relative dimension of $X_j$ over
$S$. For example, given a morphism $f:X\to Y$ in $\SP (S)$, the
transpose $\Gamma _f^{\rm t}$ of its graph $\Gamma _f$ is in $\Corr
_S^0(X,Y)$. For any two correspondences $f:X\to Y$ and $g:Y\to Z$
their composition $g\circ f$ is defined, as usual, by the formula
  $$
  g\circ f={p_{13}}_*(p_{12}^*(f)\cdot p_{23}^*(g))\; ,
  $$
where the central dot denotes the intersection of cycle classes in
the sense of \cite{Ful}. The category $\cm (S)$ of Chow-motives over
$S$ with coefficients in $\QQ $ can be defined then as a
pseudoabelian envelope of the category of correspondences with
certain ``Tate twists" indexed by integers (see \cite{Jannsen1}).
For any smooth projective $X$ over $S$ its motive $M(X/S)$ is
defined by the relative diagonal $\Delta _{X/S}$, and for any
morphism $f:X\to Y$ in $\SP (S)$ the correspondence $\Gamma _f^{\rm
t}$ defines a morphism $M(f):M(Y/S)\to M(X/S)$. The category $\cm
(S)$ is rigid with a tensor product satisfying the formula
  $$
  M(X/S)\otimes M(Y/S)=M(X\times _SY)\; ,
  $$
so that the functor
  $$
  M:\SP (S)^{\rm op}\lra \cm (S)
  $$
is tensor. The scheme $S/S$ indexed by $0$ gives the unite $\uno _S$
in $\cm (S)$, and when it is indexed by $-1$, it gives the Lefschetz
motive $\Le _S$. If $E$ is a multisection of degree $w>0$ of $X/S$,
we set
  $$
  \pi _0=\frac {1}{w}\, [E\times _SX]\; \; \; \hbox{and}\; \; \;
  \pi _{2e}=\frac {1}{w}\, [X\times _SE]\; ,
  $$
where $e$ is the relative dimension of $X/S$. Then one has the
standard isomorphisms $\uno _S\cong (X,\pi _0,0)$ and $\Le
_S^{\otimes e}\cong (X,\pi _{2e},0)$. Finally, if $f:T\to S$ is a
morphism of base schemes over $k$, then $f$ gives a base change
tensor functor $f^*:\cm (S)\to \cm (T)$. All the details about Chow
motives over a smooth base can be found, for instance, in \cite{GM}.

\subsection{\it Finite dimensional objects}
\label{findim}

Below we will use basic facts from the theory of finite dimensional
motives, or, more generally, finite dimensional objects, see
\cite{Kimura} or \cite{AK}. Roughly speaking, as soon as we have a
tensor pseudoabelian category $\bcC $ with coefficients in a field
of characteristic zero, one may speak about wedge and symmetric
powers of any object in $\bcC $. Then we say that $X\in \Ob (\bcC )$
is finite dimensional, \cite{Kimura}, if it can be decomposed into a
direct sum, $X=Y\oplus Z$, such that $\wedge ^mY=0$ and $\Sym ^nZ=0$
for some non-negative integers $m$ and $n$. The property to be
finite dimensional is closed under direct sums and tensor products,
etc. A morphism $f:X\to Y$ in $\bcC $ is said to be numerically
trivial if for any morphism $g:Y\to X$ the trace of the composition
$g\circ f$ is equal to zero, \cite[2.3]{Andre}. The important role
in the below arguments is played by the following result:

\begin{proposition}
\label{AKnilp} Let $\bcC $ be a tensor pseudoabelien category with
coefficients in $\QQ $. Let $X$ be a finite dimensional object in
$\bcC $ and let $f$ be a numerically trivial endomorphism of $X$.
Then $f$ is nilpotent in the ring of endomorphisms of $X$.
\end{proposition}

\begin{pf}
See \cite[7.5]{Kimura} for Chow motives and \cite[9.1.14]{AK} in the
abstract setting.
\end{pf}

\begin{lemma}
\label{traces} Let $\Xi :\bcC _1\to \bcC _2$ be a tensor functor
between two rigid tensor pseudoabelien categories, both with
coefficients in $\QQ $. Assume $\Xi $ induces an injection $\End
(\uno _{\bcC _1})\hookrightarrow \End (\uno _{\bcC _2})$. Then, if
$X$ is a finite dimensional object in $\bcC _1$ and $\Xi (X)=0$, it
follows that $X=0$ as well.
\end{lemma}

\begin{pf}
Let $g$ be an endomorphism of $X$. Then $F(\tr (g\circ 1_X))=\tr
(F(g\circ 1_X))$, see \cite{DeligneMilne}, page 116. Since $F(X)=0$,
we have $F(g\circ 1_X)=0$. Then $\tr (g\circ 1_X)=0$ because $F$ is
an injection on the rings of endomorphisms of units. Hence, the
identity morphism $1_X$ is numerically trivial. Since $X$ is finite
dimensional, it is trivial by Proposition \ref{AKnilp}.
\end{pf}

\begin{lemma}
\label{transfer} Let $f:T\to S$ be a finite \'etale cover of the
base $S$, and let $f^*:\cm (S)\to \cm (T)$ be the pull-back on the
corresponding categories of Chow motives. For any object $M\in \cm
(S)$, the Chow motive $M$ is finite dimensional if and only if the
motive $f^*M$ is finite dimensional.
\end{lemma}

\begin{pf}
The detailed proof is given in \cite{Gu}.
\end{pf}

\subsection{\it Spreads} \label{spreads}

The notion of a spread is described, for instance, in
\cite[4.2]{GreenGriffiths}. Here we recall and adapt it for our
goals. Let $\gamma :X\to S$ be a projective family over an
irreducible projective base $S$. Let $V$ be a closed irreducible
subvariety in the generic fiber $X_{\eta }$ defined over the field
$k(S)$. Since the family is projective, the generic fiber can be
embedded into a projective space $\PR ^n_{\eta }$ as a Zariski
closed subset. Then $V$ can also be considered as a closed subset in
the same $\PR ^n_{\eta }$. Let
  $$
  F_i(x)=F_i(x_0:\dots :x_n)=0
  $$
be a finite collection of homogeneous equations with coefficients in
$k(S)$ defining $V$ in $\PR ^n_{\eta }$. At the same time, the
variety $S$ being projective over $k$ can be embedded into a
projective space $\PR ^m_k$. Let $U$ be the intersection of $S$ with
an affine chart in $\PR ^m_k$. Being a closed irreducible subvariety
in $\AF ^m_k$, it can be defined by equations
  $$
  G_j(y)=G_j(y_1,\dots ,y_m)=0
  $$
with coefficients in $k$. Any coefficient of $F_i$ can be
represented then as a fraction $\frac{f}{g}$, where $f$ and $g$ are
elements in the coordinate ring $k[U]$ of the affine variety $U$.
Multiplying all the coefficients of the polynomials $F_i$ on the
product of their denominators, we can assume that they have
coefficients in $k[U]$. Let $Y=\gamma ^{-1}(U)$. Any element in
$k[U]$ can be represented by a polynomial in variables $y_1,\dots
,y_m$ with coefficients in $k$. Substituting polynomials instead of
coefficients in $F_i(x)$ we obtain the equations
  $$
  F_i(x,y)=F_i(x_0:\dots :x_n;y_1,\dots ,y_m)=0\; .
  $$
Then define a quasi-projective variety
  $$
  \spd V_U
  $$
in $\PR ^n_k\times \AF ^m_k$ by the system of equations
  $$
  \left\{
  \begin{array}{l} F_i(x;y)=0 \\
  G_j(y)=0
  \end{array}
  \right.
  $$
Certainly, it is a closed subvariety in $Y$, and, as such, it is a
quasi-projective subvariety in $X$. This $\spd V$ may be called a
spread of $V$ over $U$ with respect to the map $\gamma $. If we
choose another affine chart $W$ in $S$, then we will construct a
spread $\spd V_W$ over $W$. If we fix polarizations of varieties and
system of defining equations, the local spreads $\spd V_U$, $\spd
V_W$, etc are coherent, so that we can glue them into the global
spread
  $$
  \spd V
  $$
of $V$ w.r.t. $\gamma $. However, the construction is not defined
uniquely because of a few ambiguities -- in the choice of
polarizations of the varieties $X$ and $S$, and the equations
$F_i=0$ and $G_j=0$. If $X\to S$ is a family over a quasi-projective
base $S$, then we can consider any closure $\bar S$, build a spread
of $V$ with respect to the map $X\to \bar S$ and restrict it on the
preimage of the intersection $U\cap S$ in $\bar S$.

By linearity, one can define now a spread of any given algebraic
cycle $Z=\sum _in_iZ_i$ on the generic fiber $X_{\eta }$ provided
all the components $Z_i$ are rational over $k(S)$:
  $$
  \spd Z=\sum _i\spd Z_i\; .
  $$
Due to the above ambiguities, one algebraic cycle $Z$ can have two
spreads $\spd Z'$ and $\spd Z''$. However, the difference $\spd
Z'-\spd Z''$ is an algebraic cycle whose projection on $S$ w.r.t.
the map $\gamma $ is a proper subvariety in $S$, see
\cite{GreenGriffiths}. For example, if $X\to S$ is a family of
surfaces over a curve, then the ambiguity in spreads of a divisor on
$X_{\eta }$ is a vertical two-dimensional cycle, i.e. it is
generated by fibers of the map $\gamma $.

\section{Splitting $\uno $ and $\Le ^{\otimes e}$}
\label{unole}

Let $S$ be a smooth connected quasi-projective variety over a
field $k$, let
  $$
  \gamma :X\lra S
  $$
be a smooth projective family over $S$, let $d=\dim (S)$ and let
$e=\dim (X/S)$ be the relative dimension of $X$ over $S$. For any
point $s\in S$ let $X_s$ be the scheme-theoretical fiber of the
morphism $\gamma $. Then $X_s$ is a smooth projective variety over
the residue field $k(s)$. If $\eta $ is the generic point of $S$,
then $k(\eta )=k(S)$ and $X_{\eta }=X\times _S\eta $ is the generic
fiber of the family $\gamma $.

Assume $X_{\eta }$ has a rational point over $\eta $. Then, locally
on the base, the structural morphism $\gamma $ has a section, i.e.
there exists a Zariski open subset $U\subset S$ with $Y=\gamma
^{-1}(U)$, and a morphism
  $$
  \sigma :U\lra Y\; ,
  $$
such that $\gamma \circ \sigma =\id _U$. Let $\bar E$ be the Zariski
closure of $E=\sigma (U)$ in $X$. Since the structure map $\gamma $
is proper, $\gamma (\bar E)$ is a Zariski closed subset in $S$
containing $U$. Then $\gamma (\bar E)=S$ because the closure of $U$
is the whole $S$. Since $U$ is an open subset in the irreducible
variety $S$, $U$ is irreducible itself. Therefore, $\bar E$ is also
irreducible. Moreover, $\dim (\bar E)=\dim (E)=d$. Then $\bar E$ is
a one-dimensional cycle of degree one over $S$, and we consider two
relative projectors
  $$
  \bar \pi _0=[\bar E\times _SX]\hspace{5mm}\hbox{and}\hspace{5mm}
  \bar \pi _{2e}=[X\times _S\bar E]\; .
  $$
If $[\bar E]\cdot [\bar E]=0$, then $\bar \pi _0$ and $\bar \pi _2$
are pair-wise orthogonal. Assume that
  $$
  [\bar E]\cdot [\bar E]\neq 0\; .
  $$
In that case $\bar \pi _0$ and $\bar \pi _2$ are not orthogonal, so
that we cannot use them in order to split the unit and the Lefschetz
motives from $M(X/S)$ simultaneously. Let
  $$
  \zeta =\gamma _*([\bar E]\cdot [\bar E])\; .
  $$
For any natural $n$ let
  $$
  \gamma _n:X\times _S\dots \times _SX\lra S\; ,
  $$
where the fibered product is taken $n$-times. The pull-back
  $$
  \theta =\gamma _2^*(\zeta )
  $$
is a cycle class in $CH_{d+1}(X\times _SX)$. It can be considered as
a vertical correcting term for the projector $\bar \pi _0$ in the
following sense. Set
  $$
  \tau _0=\bar \pi _0-\theta \; .
  $$
As we will show right now, $\tau _0$ and $\bar \pi _{2e}$ are now
pairwise orthogonal, so that we can use them for simultaneous
splitting of $\uno $ and $\Le ^{\otimes e}$ .

For any $i,j\in \{ 1,2,3\} $ let
  $$
  p_{ij}:X\times _S\dots \times _SX\lra X\times _SX
  $$
and
  $$
  p_i:X\times _S\dots \times _SX\lra X
  $$
be the projections corresponding to their indexes. From the
commutative diagram
  $$
  \diagram
  X\times _SX\times _SX \ar[rrdd]^-{\gamma _3}
  \ar[rr]^-{p_{ij}} & & X\times _SX \ar[dd]^-{\gamma _2} \\ \\
  & & S
  \enddiagram
  $$
one has
  $$
  p_{ij}^*(\theta )=
  p_{ij}^*\gamma _2^*(\zeta )=\gamma _3^*(\zeta )\; .
  $$
Therefore,
  $$
  \begin{array}{rcl}
  \theta \circ \theta
  &=&
  {{p_{13}}_*}(p_{12}^*(\theta )\cdot p_{23}^*(\theta ))\\
  &=&
  {{p_{13}}_*}\gamma _3^*(\zeta \cdot \zeta )\\
  &=&
  {{p_{13}}_*}{{p_{13}}^*}\gamma _2^*(\zeta \cdot \zeta )\\
  &=&
  \gamma _2^*(\zeta \cdot \zeta )\cdot {{p_{13}}_*}([X\times _SX\times _SX])\\
  &=&
  \gamma _2^*(\zeta \cdot \zeta )\cdot 0\\
  &=&
  0\; .
  \end{array}
  $$
The commutative diagram
  $$
  \diagram
  X\times _SX\times _SX \ar[rrdd]^-{p_2}
  \ar[rr]^-{p_{23}} & & X\times _SX \ar[dd]^-{p_1} \\ \\
  & & X
  \enddiagram
  $$
gives
  $$
  p_{23}^*(\bar \pi _0 )=
  p_{23}^*p_1^*([\bar E])=p_2^*([\bar E])\; .
  $$
Then we compute:
  $$
  \begin{array}{rcl}
  \bar \pi _0\circ \theta
  &=&
  {{p_{13}}_*}(p_{12}^*(\theta )\cdot p_{23}^*(\bar \pi _0))\\
  &=&
  {{p_{13}}_*}(\gamma _3^*(\zeta )\cdot p_2^*([\bar E]))\\
  &=&
  {{p_{13}}_*}(p_2^*\gamma _1^*(\zeta )\cdot p_2^*([\bar E]))\\
  &=&
  {{p_{13}}_*}p_2^*(\gamma _1^*(\zeta )\cdot [\bar E])\\
  &=&
  \gamma _2^*{\gamma _1}_*(\gamma _1^*(\zeta )\cdot [\bar E])\\
  &=&
  \gamma _2^*(\zeta \cdot {\gamma _1}_*([\bar E]))\\
  &=&
  \gamma _2^*(\zeta \cdot [S])\\
  &=&
  \gamma _2^*(\zeta )\\
  &=&
  \theta \; .\\
  \end{array}
  $$
The transposition of the cycle class $\theta $ is $\theta $, and
the transposition of the cycle class $\bar \pi _0$ is $\bar \pi
_{2e}$. Hence the last equality yields:
  $$
  \theta \circ \bar \pi _{2e}=\theta \; .
  $$
Since $\bar \pi _{2e}=[X\times _S\bar E]$, we have:
  $$
  \begin{array}{rcl}
  \bar \pi _{2e}\circ \theta
  &=&
  {{p_{13}}_*}(p_{12}^*(\theta )\cdot p_{23}^*(\bar \pi _{2e}))\\
  &=&
  {{p_{13}}_*}(\gamma _3^*(\zeta )\cdot p_3^*([\bar E]))\\
  &=&
  {{p_{13}}_*}(p_3^*\gamma _1^*(\zeta )\cdot p_3^*([\bar E]))\\
  &=&
  {{p_{13}}_*}p_3^*(\gamma _1^*(\zeta )\cdot [\bar E])\\
  &=&
  0
  \end{array}
  $$
because ${{p_{13}}_*}p_3^*=0$ in general. Transposing the cycles
we get:
  $$
  \theta \circ \bar \pi _0=0\; .
  $$
At last,
  $$
  \begin{array}{rcl}
  \bar \pi _0\circ \bar \pi _{2e}
  &=&
  {{p_{13}}_*}(p_{12}^*(\bar \pi _{2e})\cdot p_{23}^*(\bar \pi _0))\\
  &=&
  {{p_{13}}_*}(p_{12}^*p_2^*([\bar E])\cdot p_{23}^*p_1^*([\bar E]))\\
  &=&
  {{p_{13}}_*}p_2^*([\bar E]\cdot[\bar E])\\
  &=&
  \gamma _2^*{\gamma _1}_*([\bar E]\cdot[\bar E])\\
  &=&
  \gamma _2^*(\zeta )\\
  &=&
  \theta \\
  \end{array}
  $$
and
  $$
  \begin{array}{rcl}
  \bar \pi _{2e}\circ \bar \pi _0
  &=&
  {{p_{13}}_*}(p_{12}^*(\bar \pi _0)\cdot p_{23}^*(\bar \pi _{2e}))\\
  &=&
  {{p_{13}}_*}(p_{12}^*p_1^*([\bar E])\cdot p_{23}^*p_2^*([\bar E]))\\
  &=&
  {{p_{13}}_*}(p_1^*([\bar E])\cdot p_3^*([\bar E]))\\
  &=&
  {{p_{13}}_*}([\bar E\times _SX\times _S\bar E]))\\
  &=&
  0\; .\\
  \end{array}
  $$

Now, using the obtained equalities, one can directly compute:
  $$
  \begin{array}{rcl}
  \tau _0\circ \tau _0
  &=&
  (\bar \pi _0-\theta )\circ (\bar \pi _0-\theta ) \\
  &=&
  \bar \pi _0\circ \bar \pi _0
  -\bar \pi _0\circ \theta
  -\theta \circ \bar \pi _0
  +\theta \circ \theta \\
  &=&
  \bar \pi _0 -\theta -0+0 \\
  &=&
  \tau _0\; ,\\
  \end{array}
  $$
  $$
  \begin{array}{rcl}
  \tau _0\circ \bar \pi _{2e}
  &=&
  (\bar \pi _0-\theta )\circ \bar \pi _{2e} \\
  &=&
  \bar \pi _0\circ \bar \pi _{2e}-\theta \circ \bar \pi _{2e} \\
  &=&
  \theta -\theta \\
  &=&
  0 \\
  \end{array}
  $$
and
  $$
  \begin{array}{rcl}
  \bar \pi _{2e}\circ \tau _0
  &=&
  \bar \pi _{2e} \circ (\bar \pi _0-\theta )\\
  &=&
  \bar \pi _{2e} \circ \bar \pi _0 -\bar \pi _{2e}\circ \theta \\
  &=&
  0-0 \\
  &=&
  0\; . \\
  \end{array}
  $$

So, we see that $\tau _0$ and $\bar \pi _{2e}$ are pair-wise
orthogonal idempotents in the associative ring $\Corr ^0_S(X,X)$.
Let
  $$
  \tilde \pi =\Delta _{X/S}-\tau _0-\bar \pi _2\; .
  $$
Since $\tau _0$ and $\bar \pi _2$ are pair-wise orthogonal, it
follows that all the projectors $\tau _0$, $\tilde \pi $ and $\bar
\pi _{2e}$ are pair-wise orthogonal. Let
  $$
  A=(X/S,\tau _0,0)\; ,
  $$
and
  $$
  \tilde M(X/S)=(X/S,\tilde \pi ,0)
  $$
be two motives in $\cm (S)$. The projector $\bar \pi _{2e}$
defines, of course, the $e$-th tensor power of relative Lefschetz
motive $\Le _S$. Then one has:
  $$
  M(X/S)=A\oplus \tilde M(X/S)\oplus \Le _S^{\otimes e}\; .
  $$

Certainly, we can split the motive $M(X/S)$ in another way. Set
  $$
  \tau _{2e}=\bar \pi _{2e}-\theta \; ,
  $$
  $$
  \hat \pi =\Delta _{X/S}-\bar \pi _0-\tau _{2e}\; .
  $$
By analogous arguments, we have that $\bar \pi _0$, $\hat \pi $
and $\tau _{2e}$ are pair-wise orthogonal idempotents, so that we
can consider the motives
  $$
  B=(X/S,\tau _{2e},0)
  $$
and
  $$
  \hat M(X/S)=(X/S,\hat \pi ,0)
  $$
as submotives in $M(X/S)$. Then $M(X/S)$ can be decomposed as
  $$
  M(X)=\uno _S\oplus \hat M(X/S)\oplus B\; .
  $$

Finally, it is easy to show that $A$ is isomorphic to $\uno _S$
and $B$ is isomorphic to the $e$-th tensor power Lefschetz motive
$\Le $. Let us prove this assertion, for example, the case of the
motive $B$. Since $\Le ^{\otimes e}=(S/S,\Delta _{S/S},-e)$, we
consider two correspondences
  $$
  a:=[X]\in \Corr ^{-e}_S(X,S)=CH^0(X)
  $$
and
  $$
  b:=[\bar E]-\sum _{i=1}^mn_i[V_i]\in
  \Corr ^e_S(S,X)=CH^e(X)\; .
  $$
Then it is easy to compute:
  $$
  b\circ a={p_{13}}_*(p_{12}^*(a)\cdot p_{23}^*(b))=
  [X\times _S\bar E]-\sum _{i=1}^mn_i[V_i\times _{Z_i}V_i]=\tau _2
  $$
and
  $$
  a\circ b={p_{13}}_*(p_{12}^*(b)\cdot p_{23}^*(a))
  ={\gamma _1}_*\left(b\right)=\Delta _{S/S}\; .
  $$

As a result, one has the following global decomposition
  $$
  M(X/S)=\uno _S\oplus M^{\dag }(X/S)\oplus \Le _S^{\otimes e}\; ,
  $$
where $M^{\dag }(X/S)$ is either $\tilde M$ or $\hat M$.

Note that in two cases the above construction gives the complete
relative Murre decomposition of the motive $M(X/S)$. For example,
this is so if $X/S$ is a smooth projective family of curves over
$S$, or if $X/S$ is a family of surfaces over $S$, with $H^1(X_{\eta
})=H^2_{\tr }(X_{\eta })=0$. Let now $\gamma :X\to S$ be a smooth
projective family of surfaces over a smooth quasi-projective curve
$S$, and let
  $$
  \Delta _{\eta }=\pi _0+\pi _1+\pi _2+\pi _3+\pi _4
  $$
be the Chow-K\"unneth decomposition of the diagonal for the surface
$X_{\eta }$, see \cite{Murre}. Recall that all the summands are
pairwise orthogonal projectors on $X_{\eta }$, and the cycle class
homomorphism maps $\pi _i$ into $(i,4-i)$-component of the diagonal
in the K\"unneth decomposition. For each index $i$ let
  $$
  \spd \pi _i
  $$
be a spread of the projector $\pi _i$. The natural question is now
as follows: is it possible to add vertical correcting terms to the
spreads $\spd \pi _i$ getting a nice decomposition of the diagonal
$\Delta _{X/S}$. Here the word ``nice" means that we wish to get a
sum of pairwise orthogonal relative projectors, which could be
considered as a relative Chow-K\"unneth decomposition under any
reasonable cohomology of $X/S$.

\section{The proof of Theorem \ref{th1}}
\label{surfaces}

Now we are ready to prove Theorem \ref{th1}. Let $\gamma :X\to S$ be
a smooth projective one-parameter family of surfaces over a smooth
irreducible curve $S$ with $H^1(X_{\eta })=H^2_{\tr }(X_{\eta })=0$.
We first consider the case when the generic fiber $X_{\eta }$ has a
rational point over $k(S)$. In that case we have the projector
$\tilde \pi $ constructed in the previous section. Let us now denote
this projector as
  $$
  \tilde \pi _2\; ,
  $$
since it is a relative second Murre projector for the whole family
$\gamma $. Certainly, this is not a unique choice of Murre's
projector, because we can also take $\hat \pi $, which can be
denoted as
  $$
  \hat \pi _2\; .
  $$
Let then $M^2(X/S)$ be the motive defined by either $\tilde \pi _2$
or $\hat \pi _2$. As we know from the previous section, one has the
decomposition
  $$
  M(X/S)=\uno _S\oplus M^2(X/S)\oplus \Le _S^{\otimes 2}\; .
  $$

Assume now $M(X/S)$ is finite dimensional. It follows that
$M^2(X/S)$ is finite dimensional as well. But, actually, $M^2$ is
evenly finite dimensional of dimension $b_2$, where $b_2$ is the
second Betti number for the generic fiber. Indeed, apriori one has
  $$
  M^2(X/S)=K\oplus L\; ,
  $$
where $\wedge ^mK=0$ and $\Sym ^nL=0$. Consider the base change
functor
  $$
  \Xi :\cm (S)\lra \cm (\eta )\; .
  $$
The functor $\Xi $ is tensor, so that it respects finite
dimensionality. In particular, the motive of the generic fiber
  $$
  M(X_{\eta })=\Xi (M(X/S))
  $$
is finite-dimensional because the relative motive $M(X/S)$ is so. In
addition, $\Xi $ induces an isomorphism between $\End (\uno _S)=\QQ
$ and $\End (\uno _{\eta })=\QQ $. Since
  $$
  \Xi (M^2(X/S))=M^2(X_{\eta })
  $$
it follows that
  $$
  \Xi (L)=0
  $$
because $M^2(X_{\eta })$ is evenly finite dimensional. Then
  $$
  L=0
  $$
by Lemma \ref{traces}. Thus, $M^2(X/S)$ is evenly finite
dimensional. Using the the same arguments one can also show that it
can be annihilated by $\wedge ^{b_2+1}$.

Let $b_2=\dim H^2(X_{\eta })$ and let $D_1,\dots ,D_{b_2}$ be
divisors on $X_{\eta }$ generating $H^2(X_{\eta })$. Since the
motive $M(X_{\eta })$ is finite dimensional, and $X_{\eta }$ is a
surface with $H^1(X_{\eta })=H^2_{\tr }(X_{\eta })=0$, the second
piece $M^2(X_{\eta })$ in the Murre decomposition of $M(X_{\eta })$
can be computed as follows:
  $$
  M^2(X_{\eta })=\Le _{\eta }^{\oplus b_2}\; ,
  $$
see \cite[Theorem 2.14]{GP1}. Actually, these $b_2$ copies of $\Le
_{\eta }$ are arising from the collection of divisor classes
  $$
  [D_1],\dots ,[D_{b_2}]
  $$
on $X_{\eta }$, and their Poincar\'e dual
  $$
  [D'_1],\dots ,[D'_{b_2}]\; ,
  $$
loc.cit. In other words, if
  $$
  \nu _2=\Xi (\tilde \pi _2)=(\tilde \pi _2)_{\eta }
  $$
is a projector determining the middle motive $M^2(X_{\eta })$, the
difference
  $$
  \xi _{\eta }=\nu _2-\sum _{i=1}^{b_2}[D_i\times _{\eta }D_i']
  $$
is homologically trivial. Then
  $$
  \xi _{\eta }^n=0
  $$
in the associative ring $\End (M^2_{\eta })$ for some $n$ by
Kimura's nilpotency theorem.

Now let
  $$
  W_i=\spd D_i
  $$
and
  $$
  W'_i=\spd D_i'
  $$
be spreads of the above divisors $D_i$ and $D_i'$ over $S$. They are
defined not uniquely, of course. The cycles $W_i\times _SW_i'$ are
in $\Corr _S^0(X\times _SX)$, and we set
  $$
  \xi =\tilde \pi _2-
  \sum _{i=1}^{b_2}[W_i\times _SW_i']\; .
  $$
Let $\omega $ be any endomorphism of the motive $M^2(X/S)$.
  $$
  \Xi (\tr (\omega \circ \xi ))=\tr (\Xi (\omega \circ \xi ))=
  \tr (\omega _{\eta }\circ \xi _{\eta })=0
  $$
because $\xi _{\eta }$ is homologically trivial (here we use the
formula on page 116 in \cite{DeligneMilne} again). Since the functor
$\Xi $ induces an isomorphism $\End (\uno _S)\cong \End (\uno _{\eta
})$,
  $$
  \tr (\omega \circ \xi )=0
  $$
for any $\omega $, i.e. $\xi $ is numerically trivial. Therefore,
  $$
  \xi ^n=0
  $$
in $\End (M^2(X/S))$ by Proposition \ref{AKnilp}.

For any cycle class $z\in CH^i(X)$ and for any correspondence $c\in
\Corr _S^j(X,X)$ let, as usual,
  $$
  c_*(z)={p_2}_*(p_1^*(z)\cdot c)\in CH^{i+j}(X)
  $$
be the action of the correspondence $c$ on $z$. In particular, if
$z\in CH^2(X)$, one has a decomposition
  $$
  \begin{array}{rcl}
  z
  &=&
  {\Delta _{X/S}}_*(z) \\
  &=&
  (\bar \pi _0)_*(z)+(\tilde \pi _2)_*(z)+(\tau _4)_*(z) \\
  &=&
  (\bar \pi _0)_*(z)+(\tilde \pi _2)_*(z)+(\bar \pi _4)_*(z)+\theta _*(z)
  \end{array}
  $$
in $CH^2(X)$. Let us compute each term separately:
  $$
  \begin{array}{rcl}
  (\bar \pi _0)_*(z)
  &=&
  {p_2}_*(p_1^*(z)\cdot \bar \pi _0) \\
  &=&
  {p_2}_*(p_1^*(z)\cdot p_1^*[\bar E]) \\
  &=&
  {p_2}_*p_1^*(z\cdot [\bar E]) \\
  &=&
  \gamma ^*\gamma _*(z\cdot [\bar E])
  \; ,
  \end{array}
  $$
  $$
  \begin{array}{rcl}
  (\bar \pi _4)_*(z)
  &=&
  {p_2}_*(p_1^*(z)\cdot \bar \pi _4) \\
  &=&
  {p_2}_*((z\times _S[X])\cdot ([X]\times _S[\bar E])) \\
  &=&
  {p_2}_*(z\times _S[\bar E]) \\
  &=&
  \gamma _*(z)\times _S[\bar E]
  \end{array}
  $$
and
  $$
  \begin{array}{rcl}
  \theta _*(z)
  &=&
  {p_2}_*(p_1^*(z)\cdot \theta ) \\
  &=&
  {p_2}_*(p_1^*(z)\cdot \gamma _2^*(\zeta )) \\
  &=&
  {p_2}_*(p_1^*(z)\cdot p_1^*\gamma ^*(\zeta )) \\
  &=&
  {p_2}_*p_1^*(z\cdot \gamma ^*(\zeta )) \\
  &=&
  \gamma ^*\gamma _*(z\cdot \gamma ^*(\zeta )) \\
  &=&
  \gamma ^*(\gamma _*(z)\cdot \zeta )
  \; .
  \end{array}
  $$

Now assume that $z\in CH^2(X)_0$, i.e. $\gamma _*(z)=0$. In that
case:
  $$
  (\bar \pi _4)_*(z)=0
  $$
and
  $$
  \theta _*(z)=0\; .
  $$
In addition, $z\cdot [\bar E]=0$ because both cycle classes have
codimension two in a smooth threefold, so that $(\bar \pi
_0)_*(z)=0$ as well. Then,
  $$
  z=(\tilde \pi _2)_*(z)\; .
  $$
On the other hand, we know that
  $$
  \tilde \pi _2=\sum _{i=1}^{b_2}[W_i\times _SW_i']+
  \xi \; ,
  $$
whence we get:
  $$
  (\tilde \pi _2)_*(z)=
  \sum _{i=1}^{b_2}([W_i\times _SW_i']_*(z))
  +\xi _*(z)\; .
  $$
Let
  $$
  v_1=-\sum _{i=1}^{b_2}([W_i\times _SW_i']_*(z))
  $$
and write $z$ as
  $$
  z=\left[\sum _jn_jZ_j\right]\; ,
  $$
where the cycles $Z_j$ are irreducible curves on $X$. For any $i$
and $j$ one has:
  $$
  [W_i\times _SW_i']_*[Z_j]=
  {p_2}_*([Z_j\times _SX]\cdot [W_i\times _SW_i'])\; .
  $$
Since $Z_j$ is of codimension two and $W_i$ is of codimension one in
the threefold $X$, it follows that a most intersection of $Z_j\times
_SX$ with $W_i\times _SW_i'$ can be represented by a sum of
algebraic cycles of type
   $$
   p\times _{k(q)}(W_i'\cdot X_q)\; ,
   $$
where $q$ is a closed point on $S$, $X_q$ is a cycle-theoretic fiber
of the morphism $\gamma $ at $q$, and $p$ is a zero-dimensional
point in the fiber $X_q$. Hence, $[W_i\times _SW_i']_*[Z_j]$ is a
linear combination of vertical cycles of type
  $$
  W_i'\cdot X_q\; .
  $$
But then so is $v_1$. In particular, $v_1$ is a vertical cycle class
and
  $$
  \xi _*(z)=z+v_1\; .
  $$
Applying $\xi $ once again we see that
  $$
  \xi ^2_*(z)=\xi _*(z+v_1)=z+v_1+\xi _*(v_1)\; .
  $$

Now we need to show that $\xi _*(v_1)$ can be also represented by a
linear combination of cycles $W_i'\cdot X_q$. By construction,
  $$
  \xi =
  \Delta _{X/S}-\bar \pi _0+\theta -\bar \pi _2-\sum _i[W_i\times _SW_i']\; .
  $$
The one-dimensional cycle class $(\bar \pi _0)_*[W_i'\cdot X_q]$ is
represented by the push-forward of the algebraic cycle
  $$
  (\bar E\cdot W_i'\cdot X_q)\times _{k(q)}X_q
  $$
on $X\times _SX$ with respect to the projection $p_2:X\times _SX\to
X$. But the intersection
  $$
  \bar E\cdot W_i'\cdot X_q
  $$
is zero because of codimensional reasons, whence
  $$
  (\bar \pi _0)_*[W_i'\cdot X_q]=0\; .
  $$
Also it is easy to show that
  $$
  (\bar \pi _2)_*[W_i'\cdot X_q]=0
  $$
and
  $$
  \theta _*[W_i'\cdot X_q]=0\; .
  $$
The diagonal does not change anything. Finally, since we have seen
already that for any codimension two subvariety $Z$ in $X$ the cycle
class
  $$
  [W_i\times _SW_i']_*[Z]
  $$
is a linear combination of vertical cycles of type $W_i'\cdot X_q$,
it is so for $W_i'\cdot X_q$ itself. Thus, the cycle class $\xi
_*(v_1)$ can be represented by a linear combination of cycles
$W_i'\cdot X_q$. In particular, $\xi _*(v_1)$ is vertical again.

Applying $\xi $ once again we see that
  $$
  \xi ^2_*(z)=z+v_2\; ,
  $$
where
  $$
  v_2=v_1+\xi _*(v_1)
  $$
is a vertical class. And so forth. After $n$ steps we will get:
  $$
  \xi ^n_*(z)=z+v_n\; ,
  $$
where $v_n$ is represented by a linear combination of the cycles
$W_i'\cdot X_q$, and, as such, it is a vertical class cycle again.
But we know that $\xi $ is a nilpotent correspondence in $\Corr
^0_S(X,X)$, so that
  $$
  \xi ^n_*(z)=0
  $$
for big enough $n$. It follows that
  $$
  z=-v_n
  $$
is a vertical cycle class.

If the generic fiber $X_{\eta }$ has no rational points over $k(S)$,
we consider a finite normal extension $L$ of $k(S)$, such that
$X_{\eta }(L)\neq \emptyset $. Let $b:S'\to S$ be a finite \'etale
cover of the curve $S$ with $k(S')=L$, and let
  $$
  \diagram
  X' \ar[dd]_-{a} \ar[rr]^-{\gamma '} & & S' \ar[dd]^-{b} \\ \\
  X \ar[rr]^-{\gamma } & & S
  \enddiagram
  $$
be a Cartesian square. Since $b$ is flat and $\gamma $ is proper,
one has $b^*\gamma _*=\gamma '_*a^*$ (see \cite[1.7]{Ful}). Then we
have the following commutative diagram
  $$
  \xymatrix{
  \ch ^2(X')_0 \ar[rrr]^-{} \ar@<-0.5ex>[dd]_-{a_*} & & &
  \ch ^2(X') \ar[rrr]^-{\gamma '_*} \ar@<-0.5ex>[dd]_-{a_*} & & &
  \ch ^0(S') \ar@<-0.5ex>[dd]_-{b_*} \\ \\
  \ch ^2(X)_0 \ar@<-0.5ex>[uu]_-{a^*} \ar[rrr]^-{} & & &
  \ch ^2(X) \ar@<-0.5ex>[uu]_-{a^*} \ar[rrr]^-{\gamma _*} & & &
  \ch ^0(S) \ar@<-0.5ex>[uu]_-{b^*}
  }
  $$
By Lemma \ref{transfer} the motive $M(X'/S')$ is finite dimensional
because $M(X/S)$ is so. Since the generic fiber of the family
$\gamma '$ has a rational point over $k(S')$, we can apply the above
arguments to show that any cycle class in $CH^2(X')_0$ is vertical.
The morphism $a:X'\to X$ is a finite \'etale cover, whence the
composition $a_*a^*$, considered on the group $CH^2(X)_0$, is just
the multiplication by $\deg (X'/X)$. Then it is easy to see that any
cycle class in $CH^2(X)_0$ is vertical as well. So, the first part
of Theorem \ref{th1} is done.

As to the second part, assume that $S=\PR ^1$, whence $\Pic
^0(S)=0$, i.e. any two closed points $q$ and $q'$ are proportional
modulo rational equivalence on $\PR ^1$. It follows that the
corresponding fibers $X_q$ and $X_{q'}$ are proportional as divisors
in $X$. Assume $X_{\eta }(k(S))\neq \emptyset $. Then all the above
vertical cycles are linear combination of cycles $W_i'\cdot X_q$ for
some fixed closed point $q$ on $\PR ^1$. Therefore, $CH^2(X)_0=\QQ
^n$ and $n\leq b_2$, where $b_2$ is the second Betti number of the
generic fiber.

If $X_{\eta }$ has no $k(\PR ^1)$-rational points, then we consider
a finite \'etale covering $S'\to S=\PR ^1$, such that, if $\eta
'=\Spec (k(S'))$, the generic fiber $X_{\eta '}$ has a rational
point over $k(S')$.

Note that any finite \'etale covering can be dominated by Galois
one. Therefore we may even assume that $S'/S$ is Galois. Then $X'/X$
is Galois as well.

Since $k(S')$ is a finite extension of $k(S)$, we can take an
algebraic closure $\overline {k(S)}$ of $k(S)$ containing $k(S')$,
whence $\bar \eta '=\bar \eta $. Fix an embedding of $k(S')$ into
the field $\CC $. Without loss of generality, we may also assume
that the algebraic closure $\overline {k(S)}=\overline {k(S')}$ is
contained in $\CC $. Then we consider also the surface $X_{\eta
}\otimes \CC $ over $\CC $ and the corresponding commutative
diagram:
    $$
    \xymatrix{
    CH^1(X_{\eta }\otimes \CC ) \ar[rr]^-{cl_{\CC }} & &
    H^2(X_{\eta }\otimes \CC ) \\ \\
    CH^1(X_{\eta '}) \ar[rr]^-{cl'} \ar[uu]^-{} & &
    H^2(X_{\eta }) \ar[uu]_-{} \\ \\
    CH^1(X_{\eta }) \ar[rr]^-{cl} \ar[uu]^-{a^*} & &
    H^2(X_{\eta }) \ar[uu]_-{=} \\ \\
    }
    $$
Since $cl$ is a surjection, $cl'$ is a surjection. But then $cl_{\CC
}$ is a surjection as well. Since $H^1(X_{\eta })=H^1(X_{\eta
'})=0$, it follows that the Albanese variety is trivial for $X_{\bar
\eta }$. Then it it trivial also for the complex surface $X_{\eta
}\otimes \CC $. Thus, we see that $X_{\eta }\otimes \CC $ is a
complex surface with $p_g=q=0$. In that case, using the
$(1,1)$-Lefschetz theorem and the exponential sequence, one can see
that $cl_{\CC }$ is an isomorphism. Then the cycle class maps $cl$
and $cl'$ are injections. Since they are also surjective, we claim
that $cl$ and $cl'$ are isomorphisms. Then the pull-back $a^*$ is an
isomorphism as well.

In other words, we may use the divisors $a^*D_1,\dots ,a^*D_{b_2}$
as an algebraic basis for $H^2(X_{\eta '})$, where $D_i's$ are the
above algebraic basis for $H^2(X_{\eta })$. Consider the spreads
$W_i$ of $D_i$. It is easy also to see that $a^*W_i$ is a spread of
$a^*D_i$ for all $i$.


Let $z$ be any algebraic cycle in $CH^2(X)_0$. By the proven part of
Theorem \ref{th1} we have that $a^*(z)$ can be represented by a
linear combination
  $$
  \sum _im_i(a^*W_i\cdot X'_{p_i})
  $$
of cycles $a^*W_i\cdot X'_{p_i}$ where $p_i$ is a closed point on
$S'$ and $X'_{p_i}$ is a cycle-theoretic preimage of $p_i$. By the
projection formula:
  $$
  \begin{array}{rcl}
  a^*a_*a^*(z)
  &=&
  \sum _im_ia^*a_*(a^*W_i\cdot X'_{p_i}) \\
  &=&
 \sum _im_ia^*(W_i\cdot a_*(X'_{p_i})) \\
  &=&
  a^*\left(\sum _im_i(W_i\cdot a_*(X'_{p_i}))\right) \; .
  \end{array}
  $$
Since $a^*$ is injective, it follows that
  $$
  z=\frac{1}{d}\cdot a_*a^*(z)=
  \frac{1}{d}\cdot \sum _im_i(W_i\cdot a_*(X'_{p_i}))\; .
  $$
But $a_*(X'_{p_i})$ is the cycle-theoretic fiber $X_{q_i}$ of the
map $\gamma $ over the image $q_i$ of the point $p_i$ under the map
$S'\to S$ by \cite[1.7]{Ful}.

So, we see that, again, any cycle class in $CH^2(X)_0$ can be
represented by a linear combination of cycles of type $W_i\cdot X_q$
where $X_q$ is a cycle-theoretic fiber over a closed point $q\in S$.
Since $S=\PR ^1$, any two fibers $X_q$ and $X_{q'}$ are proportional
divisors on $X$, whence we can use only one fiber. This finishes the
proof of Theorem \ref{th1}.

\bigskip

{\it An application of Theorem \ref{th1}.} Note that, joint with the
localization sequence for Chow groups, Theorem \ref{th1} allows also
to compute the second Chow group in a more general situation.
Indeed, let $X$ be a smooth quasi-projective threefold defined over
a field $k$, $char(k)=0$, and assume we are given with a regular
projective flat and dominant map $\gamma :X\to S$ onto a smooth
connected quasi-projective curve $S$. We can take, for example, a
Lefschetz pencil
  $$
  \gamma :X'\lra \PR ^1
  $$
under some embedding $X\subset \PR ^n$. To compute $CH^2(X')$ is the
same as to compute $CH^2(X)$ because $X'$ is just a blow up of $X$.
Now the morphism $\gamma $ is not necessarily smooth. Suppose,
moreover, that $H^1(X_{\eta })=H^2_{\tr }(X_{\eta })=0$ and that the
motive of the generic fiber $M(X_{\eta })$ is finite dimensional.
The category of Chow motives over the generic fiber can be
considered as a colimit of categories of Chow motives over Zariski
open subsets $U$ in the base curve $S$:
  $$
  \cm (\eta )=\colim _{U\subset S}\cm (U)\; ,
  $$
where the canonical morphisms $\cm (U)\to \cm (\eta )$ are
pull-backs with respect to the morphisms $\eta \to U$. Since the
motive $M(X_{\eta })$ is finite dimensional, there exists a Zariski
subset $U$, such that, if
  $$
  Y=\gamma ^{-1}(U)\; ,
  $$
the restriction
  $$
  \gamma \vert _Y:Y\lra U
  $$
is a smooth projective family and the motive $M(Y/U)$ is finite
dimensional. In other words, we ``spread out" the finite
dimensionality of the motive $M(X_{\eta })$ over a Zariski open
subset in $S$. By Theorem \ref{th1}, the Chow group $CH^2(Y)$ is
generated by a multisection and vertical cycles. Let
  $$
  Z=X'-Y\; ,
  $$
and let
  $$
  i:Y\hra X'\; \; \; \; \; \hbox{and}\; \; \; \; \; j:Z\hra X'
  $$
be, respectively, the open and closed embeddings. Using the exact
localization sequence
  $$
  CH_1(Z)\stackrel{j_*}{\lra }CH^2(X')\stackrel{i^*}{\lra
  }CH^2(Y)\to 0
  $$
it is easy now to show that any codimension two cycle class $z$ on
$X'$ is represented by a linear combination of a multisection,
vertical cycles of type $W_i\cdot F$, where $F$ is a smooth fiber
over a closed point on $S$, and vertical one-dimensional cycles
lying in singular fibers of the map $\gamma $. In particular, if
$S=\PR ^1$ and the singular fibers have finitely generated groups of
divisors, then $CH^2(X')$ is finitely generated.

\section{Families of curves.}
\label{curves}

The relative middle projector for a smooth projective curve over a
quasi-projective base can be obtained by splitting the relative unit
and Lefschetz motives simultaneously using vertical correcting terms
$\theta $, see Section \ref{unole}. In addition, these $\theta $
gives the equality $\bar \pi _0\cdot \tilde \pi _1=0$. These two
things allow to generalize Kimura's theorem, \cite[4.4]{Kimura}, to
a smooth projective curve over an arbitrary smooth quasi-projective
base:

\begin{proposition}
\label{relcurves} Let $S$ be a smooth quasi-projective variety over
a field, let $X\to S$ be a smooth projective family of curves over
$S$, and let $g$ be the genus of its generic fiber. Then the motive
$M(X/S)$ is finite dimensional. To be more precise, it splits in
$\cm (S)$ as usual: $M(X/S)=\uno _S\oplus M^1(X/S)\oplus \Le _S$,
and
  $$
  \Sym ^{2g+1}M^1(X/S)=0\; .
  $$
\end{proposition}

\begin{pf} Let $X/S$ be a smooth projective family of relative dimension
one. Assume first that the structure morphism $\gamma :X\to S$ has a
section $\sigma :S\to X$. For any natural $n$ let
  $$
  h:\Sym ^{n-1}(X/S)\times _SX\lra \Sym ^n(X/S)
  $$
be the evident map\footnote{all symmetric products here and below
are over the base $S$}, and let
  $$
  i:\Sym ^{n-1}(X/S)\lra \Sym ^{n-1}(X/S)\times _SX
  $$
be the embedding induced by the section $\sigma $. Then let $M_n$
be the codimension one subvariety in $\Sym ^n(X/S)$, which is the
image of the composition $h\circ i$. Let also
  $$
  \mathcal M_n
  $$
be the invertible sheaf corresponding to the divisor $M_n$
respectively.

\begin{lemma}
\label{MattukScwarz} Under the above assumptions, there exists a
locally free sheaf $\mathcal E$ on the relative Jacobian
  $$
  J=\Pic ^0_{X/S}
  $$
of $X/S$, and an isomorphism
  $$
  r_n:\Sym ^n(X/S)\lra  \PR (\mathcal E)\; ,
  $$
such that the diagram
  $$
  \diagram
  \Sym ^n(X/S) \ar[ddrr]^-{} \ar[rr]^-{r_n} & &
  \PR (\mathcal E) \ar[dd]^-{} \\ \\
  & & J
  \enddiagram
  $$
commutes and
  $$
  r_n^*\mathcal O(1)=\mathcal M_n\; ,
  $$
where all three maps in the above commutative diagram are over
$S$.
\end{lemma}

\begin{pf}
Follow \cite{Schwarz} in the relative setting.
\end{pf}

If the generic fiber $X_{\eta }$ has no $k(S)$-rational points,
then, as in the proof of Theorem \ref{th1}, we take a finite normal
extension $L$ of the field $k(S)$ with $X_{\eta }(L)\neq \emptyset
$, and apply Lemma \ref{transfer}. In other words, without loss of
generality, we may assume from the beginning that $X_{\eta }$ has a
point rational over $k(S)$. Then, locally on the base, Proposition
\ref{relcurves} is a straightforward generalization of results from
\cite{Kimura}. We first recall the local situation and then prove
the theorem globally.

Since $X_{\eta }(k(S))\neq \emptyset $, there exists a Zariski open
subset $U$ in $S$, such that $\gamma :Y\to U$ has a section $\sigma
$, where $Y=\gamma ^{-1}(U)$ (see Section \ref{unole}). We can make
$U$ smaller, so that the self-intersection $[E]\cdot [E]$ is zero in
the Chow group $CH_d(Y\times _UY)$, where $E=\sigma (S)$. In that
case the correspondences $\pi _0=[E\times _UY]$ and $\pi
_{2}=[Y\times _UE]$ are pair-wise orthogonal, so that we can
directly define the middle projector for $Y/U$ by the formula
  $$
  \pi _1=\Delta _{Y/U}-\pi _0-\pi _2\; .
  $$
Then we get the decomposition
  $$
  M(Y/U)=\uno _U\oplus M^1(Y/U)\oplus \Le _U\; ,
  $$
where
  $$
  M^1(Y/U)=(Y/U,\pi _1,0)\; .
  $$

The morphism $Y\to U$ is smooth projective with a section. By Lemma
\ref{MattukScwarz} there exists a locally free sheaf $\mathcal E$ on
the Jacobian $J$ of $Y/U$, such that $\Sym ^n(Y/U)$ is isomorphic to
$\PR (\mathcal E)$ provided $n>2g-2$. And, moreover, the pull-back
of the corresponding bundle $\mathcal O(1)$ to $\Sym ^n(Y/U)$
coincides with the class of the relative divisor $\mathcal M_n$. The
equality $[E]\cdot [E]=0$ gives
  $$
  \pi _0\cdot \pi _0=0\; ,
  $$
whence
  $$
  \pi _1\cdot \pi _0=0\; .
  $$
Therefore, one can directly generalize Kimura's arguments from
\cite[\S 4]{Kimura} to show that
  $$
  \Sym ^{2g+1}M^1(Y/U)=0\; .
  $$

\bigskip

Now we want to prove the theorem globally on the base. We consider
the projectors
  $$
  \bar \pi _0=[\bar E\times _SX]\hspace{5mm}\hbox{and}\hspace{5mm}
  \bar \pi _2=[X\times _S\bar E]\; ,
  $$
introduce the vertical correcting term $\theta $, and take the
correspondences $\tau _0$ and $\tau _2$ from Section \ref{unole}.
Then we can define
  $$
  \tilde \pi _1=\Delta _{X/S}-\tau _0-\bar \pi _2
  $$
in order to get the splitting with
  $$
  M^1(X/S)=(X/S,\tilde \pi _1,0)\; .
  $$
(see Section \ref{unole}). We need to show that $M^1(X/S)$ is
finite dimensional. Let again $g$ be the genus of the generic
fiber of the structural morphism $\gamma :X\to S$. For any natural
number $n$ let $(X/S)^n_i$ be the fibered product
  $$
  X\times _S\dots \times _S\bar E\times _S\dots \times _SX\; ,
  $$
where $\bar E$ is located on the $i$-th place. Since $X/S$ has a
relative divisor of degree one over $S$ (subvariety $\bar E$), we
can again apply Kimura's arguments to the family $X/S$. Therefore,
in order to prove that $\Sym ^{2g+1}M^1(X/S)$ vanishes we have only
to show that the intersection of the cycle
  $$
  [(X/S)^{2g+1}_i\times _S(X/S)^{2g+1}]
  $$
with
  $$
  \tilde \pi _1^{(2g+1)}
  $$
is equal to zero. But
  $$
  \begin{array}{rcl}
  \bar \pi _0\cdot \tilde \pi _1
  &=&
  \bar \pi _0\cdot \Delta _{X/S}-
  \bar \pi _0\cdot \tau _0-
  \bar \pi _0\cdot \bar \pi _2 \\
  &=&
  [\bar E\times _S\bar E]-
  \bar \pi _0\cdot (\bar \pi _0-\theta )-
  [\bar E\times _S\bar E] \\
  &=&
  \bar \pi _0\cdot \bar \pi _0-\bar \pi _0\cdot \theta \; .
  \end{array}
  $$
On the other hand,
  $$
  \begin{array}{rcl}
  \bar \pi _0\cdot \theta
  &=&
  p_1^*([\bar E])\cdot \gamma _2^*{\gamma _1}_*
  ([\bar E]\cdot [\bar E]) \\
  &=&
  p_1^*([\bar E])\cdot p_1^*\gamma _1^*{\gamma _1}_*
  ([\bar E]\cdot [\bar E]) \\
  &=&
  p_1^*([\bar E]\cdot \gamma _1^*{\gamma _1}_*
  ([\bar E]\cdot [\bar E])) \\
  &=&
  p_1^*([\bar E]\cdot {p_1}_*p_2^*([\bar E]\cdot [\bar E])) \\
  &=&
  p_1^*([\bar E]\cdot [\bar E]) \\
  &=&
  p_1^*([\bar E])\cdot p_2^*([\bar E]) \\
  &=&
  \bar \pi _0\cdot \bar \pi _0\; .
  \end{array}
  $$
Therefore,
  $$
  \begin{array}{rcl}
  \bar \pi _0\cdot \tilde \pi _1
  &=&
  \bar \pi _0\cdot \bar \pi _0-\bar \pi _0\cdot \theta \\
  &=&
  \bar \pi _0\cdot \bar \pi _0-\bar \pi _0\cdot \bar \pi _0 \\
  &=&
  0\; .
  \end{array}
  $$
Then,
  $$
  [(X/S)^{2g+1}_i\times _S(X/S)^{2g+1}]\cdot
  \tilde \pi _1^{(2g+1)}=0
  $$
because $\bar \pi _0\cdot \tilde \pi _1=0$.

Certainly, we can split the motive $M(X/S)$ in another way. Set
  $$
  \hat \pi _1=\Delta _{X/S}-\bar \pi _0-\tau _2\; .
  $$
and
  $$
  M^1(X/S)=(X/S,\hat \pi _1,0)\; .
  $$
In that case the proof is analogous.
\end{pf}

Of course, Proposition \ref{relcurves} has the following standard
corollaries. Let
  $$
  \cma (S)
  $$
be the full tensor pseudoabelian subcategory in $\cm (S) $ generated
by motives of relative curves. Then, by Proposition \ref{relcurves}
and using abstract properties of finite dimensional objects, we
obtain that all objects in $\cma (S)$ are finite dimensional.
Consequently, if $M\in \cma (S)$ and $f:M\to M$ is a numerically
trivial endomorphism of the motive $M$, then $f$ is nilpotent in the
associative ring $\End (M)$ by Proposition \ref{AKnilp}. The last
assertion can help to detect algebraic cycles in families made by
relative curves.

\section{Examples}
\label{examples}

Let us now indicate where threefolds satisfying the assumptions of
the above theorem are arising from.

Let $Y$ be a smooth projective complex surface with $p_g=q=0$.
Actually, $Y$ is defined over the minimal extension $K_0$ of the
prime field $\QQ $, such that all the coefficients of equations
defining $Y$ are in $K_0$. Clearly, $K_0$ is finitely generated over
$\QQ $. Since $q=0$ for $Y/\CC $, the surface $Y$ has trivial
irregularity also over $K_0$, whence $H^1(Y)=0$ for $Y$ considered
over $K_0$. On the other hand, it might be possible that there are
no a system of divisors on $Y/K_0$ generating $H^2(Y)$. Then let $K$
be any finitely generated extension of $K_0$, such that the variety
$Y/K$ satisfies the condition $H^2_{\tr }(Y)=0$ over $K$.

The field $K$ can be viewed as the field of rational functions over
an irreducible quasi-projective variety defined over $\QQ $. If we
take $K$ so that its transcendence degree over $\QQ $ is positive,
then we can also consider $K$ as the function field $k(S')$ of an
irreducible quasi-projective curve $S'$ defined over some extension
$k$ of $\QQ $. This $k$ can have a non-trivial transcendent degree
over $\QQ $, of course. Now we spread out $Y/k(S')$ getting a family
of surfaces $\gamma :X\to S$, where $S$ is a Zariski open subset in
the curve $S'$. The generic fiber $X_{\eta }$ of the morphism
$\gamma $ is then our surface $Y$ considered over $\eta =\Spec
(k(S))$, where $k(S)=K$.

Now assume Bloch's conjecture holds for $Y/\CC $. It should be noted
that, by now, we know only four types of surfaces with $p_g=q=0$ and
known Bloch's conjecture:

\begin{itemize}

\item[(i)]
Enriques surfaces, \cite{BKL};

\item[(ii)]
the classical Godeaux surface, \cite{InoseMizukami};

\item[(iii)]
some surfaces originated by modular groups, \cite{Barlow}, and

\item[(iv)]
finite quotients of intersections of four quadrics in $\PR ^6$,
\cite{Voisin}.

\end{itemize}

Since Bloch's conjecture holds for $Y/\CC $, the motive $M(Y/\CC )$
is finite dimensional, \cite{GP2}. Finite dimensionality of a Chow
motive is, actually, a rational triviality of its determining
projector, and the last thing can be viewed as an existence of an
algebraic cycle on a variety. Therefore, we can say that $M(Y)$ is
finite dimensional over some finitely generated extension of the
primary field $\QQ $ over which a trivializing algebraic cycle is
rational. Extending $K$ more, if necessary, we can assume that the
motive $M(Y/k(S))$ is finite dimensional. But $M(Y/k(S))$ is nothing
else than the motive $M(X_{\eta })$. Spreading out finite
dimensionality of the motive $M(X_{\eta })$ over some Zariski open
subset in $S$ we get a family whose relative motive is finite
dimensional.

So, we see that any smooth projective surface $X/\CC $ with
$p_g=q=0$ and finite dimensional $M(X)$ naturally gives rise to a
family of surfaces satisfying the assumptions of Theorem \ref{th1}
by means of the above spreading construction. In particular, we can
now describe codimension $2$ algebraic cycles on threefolds arising
from spreadings of surfaces of the above four types (i)-(iv).

Another one way to introduce such families is to consider Lefschetz
pencils of appropriate threefolds in a projective space.

\begin{small}

\end{small}

\vspace{4mm}

{\footnotesize

\noindent {\sc School of Mathematics, The Institute for Advanced
Study, Einstein Drive, Princeton, NJ 08540}


\noindent {\it E-mail}: {\tt guletski@ias.edu}

}

\end{document}